\documentclass[12pt,a4paper,fleqn]{article}
\usepackage{a4wide,amsfonts,amsmath,latexsym,amssymb,euscript,graphicx,units,mathrsfs}

\usepackage{graphicx}
\usepackage{color}
\usepackage{amssymb}
\usepackage{amssymb}
\usepackage[T1]{fontenc}
\usepackage{latexsym}
\usepackage{xypic}
\usepackage{eufrak}
\usepackage{euscript}
\usepackage{amsfonts,amsmath}
\usepackage{verbatim}
\usepackage{fancyhdr}
\usepackage[english]{babel}
\usepackage{mathrsfs}
\usepackage{units}

\newtheorem{prop}{Proposition}[section]
\newtheorem{proposition}[prop]{Proposition}

\newtheorem{theorem}[prop]{Theorem}
\newtheorem{thm}[prop]{Theorem}

\renewcommand{\geq}{\geqslant}
\def\leq{\leqslant}

\newcommand{\N}{\mathbb{N}}

\newcommand{\R}{\mathbb{R}}

\def\e{\varepsilon}
\font\tenbb=msbm10

\def\1{{\mathbf{1}}}

\font\tenbb=msbm10

\def\nN{\hbox{\tenbb N}}

\def\1{{\mathbf{1}}}
\def\0.5{{\frac{1}{2}}}

\newcommand{\qed}{\nopagebreak\hspace*{\fill}
{\vrule width6pt height6ptdepth0pt}\par}

\newcounter{rea}
\begin{document}

\begin{center}
{\Large{\bf Convergence in law implies convergence in total variation for
polynomials in independent Gaussian, Gamma or Beta random variables}}
\normalsize
\\~\\ Ivan Nourdin\footnote{supported in part by the (french) ANR grant `Malliavin, Stein and Stochastic Equations with Irregular Coefficients'
[ANR-10-BLAN-0121].} (Universit\'e de Lorraine)\\
 Guillaume Poly (Universit\'e du Luxembourg)\\
\end{center}

{\small \noindent {{\bf Abstract}:
Consider a sequence of 
polynomials of bounded degree evaluated in independent Gaussian, Gamma or Beta random variables.
We show that, if this sequence converges in law to a nonconstant distribution, then $(i)$
the limit distribution is necessarily absolutely continuous with respect to the Lebesgue measure
and $(ii)$ the convergence automatically takes place in the total variation topology.
Our proof, which relies on the Carbery-Wright inequality and makes use of a diffusive Markov operator approach, extends  the results of \cite{NP} to the Gamma and Beta cases.
\\


\noindent {\bf 2000 Mathematics Subject Classification:}

\section{Introduction and main results}

The Fortet-Mourier distance between (laws of) random variables,
defined as
\begin{equation}\label{fm}
d_{FM}(F,G)=
\sup_{\substack{\|h\|_\infty\leq 1\\ \|h'\|_\infty\leq 1}}
\big|
E[h(F)]-E[h(G)]
\big|
\end{equation}
is well-known to metrize the convergence in law,  see, e.g., \cite[Theorem 11.3.3]{Dudley book}. In other words,
as $n\to\infty$ one has that
$F_n\overset{\rm law}{\to}F_\infty$ if and only if $d_{FM}(F_n,F_\infty)\to 0$.
But there is plenty of other distances that allows one to measure
the proximity between laws of random variables.
For instance, one may use the Kolmogorov distance:
\[
d_{Kol}(F,G)=\sup_{x\in\R}\big| P(F\leq x)-P(G\leq x)\big|.
\]
Of course,  if $d_{Kol}(F_n,F_\infty)\to 0$ then $F_n\overset{\rm law}{\to}F_\infty$. 
But the converse implication is wrong in general, meaning that the Kolmogorov distance does {\it not} metrize the convergence in law. Nevertheless, it becomes true when the target law is {\it continuous}
(that is, when the law of $F_\infty$ has a density with respect to the Lebesgue measure), a fact which can be easily checked by using (for instance) the second Dini's theorem. 
Yet another popular distance for measuring the distance between laws of random variables, which is even stronger than the Kolmogorov distance, is the total variation distance:
\begin{equation}\label{tv1}
d_{TV}(F,G)=\sup_{A\in\mathcal{B}(\R)}\big| P(F\in A)-P(G\in A)\big|.
\end{equation}
One may prove that
\begin{equation}\label{tv2}
d_{TV}(F,G)=\frac12\,\,\,\sup_{\|h\|_\infty\leq 1}
\big|
E[h(F)]-E[h(G)]
\big|,
\end{equation}
or, whenever $F$ and $G$ both have a density (noted $f$ and $g$ respectively)
\[
d_{TV}(F,G)=\frac12\int_\R |f(x)-g(x)|dx.
\]
Unlike the Fortet-Mourier or Kolmogorov distances, 
it can happen that $F_n\overset{\rm law}{\to}F_\infty$ for continuous $F_n$ and $F_\infty$ without having that 
$d_{TV}(F_n,F_\infty)\to 0$. For an explicit counterexample, one may consider 
$F_n\sim \frac{2}{\pi}\cos^2(nx){\bf 1}_{[0,\pi]}(x)dx$; indeed, it is straightforward to check that $F_n\overset{\rm law}{\to}F_\infty\sim\mathcal{U}_{[0,\pi]}$
but $d_{TV}(F_n,F_\infty)\not\to 0$ (it is indeed a strictly positive quantity that does not depend on  $n$).
\\\\
As we just saw, the convergence in total variation is very strong and
therefore it cannot be expected from
the mere convergence in law without additional structure. Let us give three representative results in this direction. Firstly, there is a celebrated theorem of Ibragimov (see, e.g., Reiss \cite{reiss}) according to which, if $F_n, F_\infty$ 
are continuous random variables with densities
$f_n$,$f_\infty$ that
that are
{\it unimodal}, then
$F_n\overset{\rm law}{\to}F_\infty$ if and only if $d_{TV}(F_n,F_\infty)\to 0$. Secondly, let us quote the paper \cite{mamatov}, in which necessary and sufficient conditions are given (in term of the absolute continuity of the laws) so that the classical Central Limit Theorem holds in total variation. Finally, let us mention \cite{Viboga} or \cite{poly} for conditions ensuring the convergence in total variation for random variables in Sobolev or Dirichlet spaces. Although all the above examples are related to very different frameworks, they have in common the use of a particular structure of the involved variables;
loosely speaking, this structure allows to derive a kind of ``non-degeneracy'' in an appropriate sense which, in turn, enables to reinforce the convergence, from the Fortet-Mourier distance to the total variation one.
\\\\
Our goal in this short note is to exhibit yet another instance
when convergence in law and in total variation are equivalent. 
More precisely, we shall prove the following result,
which may be seen as an extension to the Gamma and Beta cases of our previous results in \cite{NP}.

\begin{theorem}\label{Main}
Assume that one of the following three conditions is satisfied:
\begin{enumerate}
\item[(1)] $X\sim N(0,1)$;
\item[(2)] $X\sim \Gamma(r,1)$ with $r\geq 1$;
\item[(3)] $X\sim \beta(a,b)$ with $a,b\geq 1$.
\end{enumerate}
Let $X_1,X_2,\ldots$ be independent copies of $X$.
Fix an integer $d\geq 1$ and,
for each $n$, let $m_n$ be a positive integer and let $Q_n\in \R[x_1,...,x_{m_n}]$ be a multilinear
polynomial of degree at most $d$; assume further that $m_n\to \infty$ as $n\to\infty$.
Finally, suppose that $F_n$ has the form
\[
F_n = Q_n(X_1,\ldots,X_{m_n}),\quad n\geq 1,
\]
and that it converges
in law as $n\to\infty$ to a non-constant random variable $F_\infty$.
Then the law of $F_\infty$ is absolutely continuous with respect to the Lebesgue measure and $F_n$ actually converges to $F_\infty$ in total variation.
\end{theorem}

In the statement of Theorem \ref{Main}, by `multilinear polynomial of degree at most $d$' we mean
a polynomial $Q\in\R[x_1,\ldots,x_m]$ of the form
\[
Q(x_1,\ldots,x_m)=\sum_{S\subset \{1,\ldots,m\},\,|S|\leq d}
a_S \,\,\,\prod_{i\in S}x_i,
\]
for some real coefficients $a_S$ and with the usual convention
that $\prod_{i\in \emptyset}x_i=1$.

Before providing the proof of Theorem \ref{Main}, let us
comment a little bit why we are `only' considering
the three cases (1), (2) and (3). 
This is actually due to our method of proof.
Indeed, the two main ingredients we are using for showing Theorem \ref{Main} are the following.
\begin{enumerate}
\item[(a)] We will make use of a Markov semigroup approach. More specifically, our strategy rely to the use of orthogonal polynomials which are also eigenvectors of diffusion operators. And in dimension 1, up to affine transformations only the Hermite (case (1)), Laguerre (case (2)) and Jacobi (case (3)) polynomials are of this form, see \cite{Mazet}. 
\item[(b)] We will make use of the Carbery-Wright inequality (Theorem \ref{cw-thm}). The main assumption for this inequality to hold is the log-concavity property. This impose some further (weak) restrictions on the parameters in the cases (2) and (3).
\end{enumerate}

The rest of the paper is organized as follows.
In Section 2, we gather some useful preliminary results.
Theorem \ref{Main} is shown in Section 3. 
\section{Preliminaries}

From now on,
we shall write $m$ instead of $m_n$
for the sake of simplicity.

\subsection{Markov semigroup}\label{preli}

In this section, we introduce the framework we shall need along the proof of Theorem \ref{Main}.
We refer the reader to \cite{Ledoux-bakry-gentil} for the details and missing proofs.
Fix an integer $m$ and let $\mu$ denote the distribution of the random vector $(X_1,\ldots,X_{m})$,
with $X_1,\ldots,X_m$ independent copies of $X$, 
for $X$ satisfying either (1), (2) or (3).
In these three cases, there exists a reversible Markov process on $\R^m$, with semigroup $P_t$, equilibrium measure $\mu$ and generator $\mathcal{L}$. 
The operator $\mathcal{L}$ is selfadjoint and
negative semidefinite. We define the Dirichlet form $\mathcal{E}$ associated to $\mathcal{L}$:
\[
\mathcal{E}(f,g)=-\int f\mathcal{L}g d\mu = -\int g\mathcal{L}f d\mu.
\]
When $f=g$, we simply write $\mathcal{E}(f)$ instead of
$\mathcal{E}(f,f)$.
The carr\'e du champ operator $\Gamma$ will be also of interest; it is the operator defined as
\[
\Gamma(f,g) = \frac12\big(\mathcal{L}(fg)-f\mathcal{L}g-g\mathcal{L}f\big).
\]
Similarly to $\mathcal{E}$, when $f=g$ we simply write $\Gamma(f)$ instead of
$\Gamma(f,f)$.
Since $\int \mathcal{L}f \,d\mu=0$, 
we observe the following link between the Dirichlet form $\mathcal{E}$ and the carr\'e du champ operator $\Gamma$:
\[
\int \Gamma(f,g) d\mu = \mathcal{E}(f,g).
\]
An important property which is satisfied in the three cases $(1)$, $(2)$ and $(3)$ is that $\Gamma$ is diffusive in the following sense: 
\begin{equation}\label{diffu}
\Gamma(\phi(f),g)=\phi'(f)\Gamma(f,g).
\end{equation}
Besides, and it is another important property shared by $(1)$, $(2)$, $(3)$, the eigenvalues of $-\mathcal{L}$
may be ordered as a countable sequence like 
$0 = \lambda_0 < \lambda_1 < \lambda_2 <  \cdots$, with a corresponding sequence of orthonormal eigenfunctions $u_0$, $u_1$, $u_2$, $\cdots$ where $u_0 =1$; in addition, this sequence of eigenfunctions forms a complete orthogonal basis of $L^2(\mu)$.
For completeness, let us give more details in each of our three cases (1), (2), (3) of interest.

\begin{itemize}
\item[(1)]
The case where $X\sim N(0,1)$. We have
\begin{equation}\label{O-U}
\mathcal{L}f(x)=\Delta f(x)-x\cdot \nabla f(x),\quad x\in \R^m,
\end{equation}
where $\Delta$ is the Laplacian operator and $\nabla$ is the gradient.
As a result,
\begin{equation}\label{gamma1}
\Gamma(f,g)=\nabla f \cdot \nabla g.
\end{equation}
We can compute that $\text{Sp}(-\mathcal{L})=\nN$ and that
$\text{Ker}(\mathcal{L}+k\,I)$ (with $I$ the identity operator) is composed of those polynomial
functions $R(x_1,\ldots,x_m)$ having the form
\[
R(x_1,\ldots,x_m)=\sum_{i_1+i_2+\cdots+i_{m}=k}\alpha(i_1,\cdots,i_{m})\prod_{j=1}^{m} H_{i_j}(x_j).
\]
Here, $H_i$ stands for the Hermite polynomial of degree $i$.

\item[(2)]
The case where $X\sim\Gamma(r,1)$. The density of $X$ is $f_X(t)=t^{r-1}\frac{e^{-t}}{\Gamma(r)}$, $t\geq 0$, with $\Gamma$ the Euler Gamma function; it is log-concave for $r\geq 1$. Besides, we have
\begin{equation}\label{Lag2}
\mathcal{L}f(x)=\sum_{i=1}^{m} \Big{(}x_{i} \partial_{ii}f+(r+1-x_i)\partial_i f\Big{)},\quad x\in\R^m.
\end{equation}
As a result,
\begin{equation}\label{gamma2}
\Gamma(f,g)(x)=
\sum_{i=1}^{m} x_i \partial_i f(x)\partial_i g(x), \quad x\in \R^m.
\end{equation}
We can compute that $\text{Sp}(-\mathcal{L})=\nN$ and
that
$\text{Ker}(\mathcal{L}+k\,I)$ is composed of those polynomial
functions $R(x_1,\ldots,x_m)$ having the form
\[
R(x_1,\ldots,x_m)=\sum_{i_1+i_2+\cdots+i_{m}=k}\alpha(i_1,\cdots,i_{m})\prod_{j=1}^{m} L_{i_j}(x_j).
\]
Here $L_i(X)$ stands for the $i$th Laguerre polynomial of parameter $r$ defined as
\[
L_i(x)= {x^{-r} e^x \over i!}{d^i \over dx^i} \left\{e^{-x} x^{i+r}\right\},\quad x\in\R.
\]

\item[(3)]
The case where $X\sim \beta(a,b)$. In this case,
$X$ is continuous with density
\[
f_X(t) = \begin{cases}\displaystyle \frac{t^{a-1}(1-t)^{b-1}}{\int_0^1 u^{a-1} (1-u)^{b-1}\, du} & \hbox{ if }t\in[0,1]\\0 & \hbox{ otherwise }\end{cases}.
\]
The density $f_X$ is log-concave when $a,b\geq 1$. Moreover, 
we have
\begin{equation}\label{Jac}
\mathcal{L}f(x)=\sum_{i=1}^{m} \Big{(}(1-x_{i}^2) \partial_{ii}f+(b-a-(b+a) x_i)\partial_i f\Big{)},\quad x\in\R^m.
\end{equation}
As a result,
\begin{equation}\label{gamma3}
\Gamma(f,g)(x)=
\sum_{i=1}^{m} (1-x_i^2) \partial_i f(x)\partial_i g(x), \quad x\in \R^m.
\end{equation}
Here, the structure of the spectrum turns out to be a little bit more complicated than in the two previous cases (1) and (2). Indeed, we have that
\[
\text{Sp}(-\mathcal{L})=\{i_1(i_1+a+b-1)+\cdots+i_{m}(i_{m}+a+b-1)\,|\,i_1,\ldots,i_m\in \nN\}.
\]
Note in particular that the first nonzero element of $\text{Sp}(-\mathcal{L})$
is $\lambda_1=a+b-1>0$.
Also, one can compute that, when $\lambda\in \text{Sp}(-\mathcal{L})$, then
$\text{Ker}(\mathcal{L}+\lambda\,I)$ is composed of those polynomial
functions $R(x_1,\ldots,x_m)$ having the form
\[
R(x_1,\ldots,x_m)=\sum_{i_1(i_1+a+b-1)+\cdots+i_{m}(i_{m}+a+b-1)=\lambda}\alpha(i_1,\cdots,i_{n_m}) J_{i_1}(x_1)\cdots J_{i_{m}}(x_{m}).
\]
Here $J_{i}(X)$ is the $i$th Jacobi polynomial, defined as
\[
J_{i}(x)=\frac{(-1)^i}{2^i i!} (1-x)^{1-a} (1+x)^{1-b} \frac{d^i}{dx^i} \left\{(1-x)^{a-1}(1+x)^{b-1}(1 - x^2)^i\right\},\quad x\in\R.
\] 
\end{itemize}
To end up with this quick summary, we stress that 
a Poincar\'e inequality holds true in  the three cases (1), (2) and (3). This is well-known and easy to prove, by using
the previous facts together with
the decomposition 
\[
L^2(\mu) = \bigoplus_{\lambda\in \text{Sp}(-\mathcal{L})}
\text{Ker}(\mathcal{L}+\lambda\,I).
\]
Namely, with $\lambda_1>0$ the first nonzero eigenvalue of
$-\mathcal{L}$, we have
\begin{equation}\label{poincare}
{\rm Var}_\mu(f)\leq \frac{1}{\lambda_1}\,\mathcal{E}(f).
\end{equation}

\subsection{Carbery-Wright inequality}
The proof of Theorem \ref{Main} will rely, among others, on the following crucial inequality due to Carbery and Wright (\cite[Theorem 8]{CW}). We state it here for convenience.

\begin{thm}[Carbery-Wright]\label{cw-thm}
There exists an absolute constant $c>0$ such that, if $Q:\R^m\to\R$ is a polynomial of degree at most $k$
and $\mu$ is a log-concave probability measure on $\R^m$ then,
for all $\alpha>0$, 
\begin{equation}\label{cw-ineq}
\left(\int Q^2d\mu\right)^{\frac1{2k}}\times 
\mu\{x\in\R^m:\,|Q(x)|\leq \alpha\}
\leq c\,k\,\alpha^{\frac1k}.
\end{equation}
\end{thm}

\subsection{Absolute continuity}
There is a celebrated result of Borell \cite{borell} according to which, if $X_1$, $X_2$, ... are independent, identically distributed and $X_1$ has an absolute continuous law, then
any nonconstant polynomial in the $X_i$'s has an absolute continuous law, too.
In the particular case where the common law satisfies either (1), or (2) or (3) in Theorem \ref{Main}, one can
recover Borell's theorem as a consequence of the previous Carbery-Wright inequality.
We provide the proof of this fact here, since it may be seen as a first step towards
the proof of Theorem \ref{Main}.

\begin{proposition}\label{abso}
Assume that one of the three conditions (1), (2) or (3) of Theorem \ref{Main} is satisfied.
Let $X_1,X_2,\ldots$ be independent copies of $X$.
Consider two integers $m,d\geq 1$ and let $Q\in \R[x_1,...,x_{m}]$ be a polynomial of degree $d$.
Then the law of $Q(X_1,\ldots,X_{m})$ is absolutely continuous with respect to the Lebesgue measure if and only if its variance is not zero.
\end{proposition}
{\it Proof}.
Write $\mu$ for the distribution of $(X_1,\ldots,X_m)$ and assume that ${\rm Var}(Q(X_1,\ldots,X_{m}))>0$.
We shall prove that, if $A$ is a Borel set of $\R$ with Lebesgue
measure zero, then $P(Q(X_1,\ldots,X_{m})\in A)=0$.
This will be done in three steps.\\

{\bf Step 1}. Let $\e>0$ and let $B$ be a {\it bounded} Borel set.
We shall prove that
\begin{equation}\label{cl1}
\int {\bf 1}_{\{Q\in B\}} 
\frac{\Gamma(Q)}{\e+\Gamma(Q)}
d\mu
=
\int \left(
\int_{-\infty}^Q {\bf 1}_B(u)du\times
\left\{
\frac{-\mathcal{L}Q}{\Gamma(Q)+\e}+2\frac{\Gamma(\Gamma(Q))}{(\Gamma(Q)+\e)^2}
\right\}
\right)d\mu.
\end{equation}
Indeed, let $h:\R\to [0,1]$ be $C^\infty$ with compact support. We can write, using among other (\ref{diffu}),
\begin{eqnarray*}
&&\int \left(
\int_{-\infty}^Q h(u)du\times
\frac{-\mathcal{L}Q}{\Gamma(Q)+\e}
\right)d\mu=
\mathcal{E}\left(
\int_{-\infty}^Q h(u)du\times
\frac{1}{\Gamma(Q)+\e},Q
\right)\\
&=&\int
\left(
h(Q)\frac{\Gamma(Q)}{\Gamma(Q)+\e}-
2\int_{-\infty}^Q h(u)du\,
\frac{\Gamma(\Gamma(Q))}{(\Gamma(Q)+\e)^2}
\right)
d\mu.
\end{eqnarray*}
Applying Lusin's theorem allows one, by dominated convergence, to pass from $h$ to
${\bf 1}_B$ in the previous identity; this leads to the desired conclusion (\ref{cl1}).\\

{\bf Step 2}. Let us apply (\ref{cl1}) to $B=A\cap [-n,n]$. Since
$\int_{-\infty}^\cdot {\bf 1}_B(u)du$ is zero almost everywhere, one
deduces that, for all $\e>0$ and all $n\in\N^*$,
\[
\int {\bf 1}_{\{Q\in A\cap [-n,n]\}} 
\frac{\Gamma(Q)}{\e+\Gamma(Q)}
d\mu=0.
\] 
By monotone convergence $(n\to\infty)$ it comes that,  for all $\e>0$,
\begin{equation}\label{bla}
\int {\bf 1}_{\{Q\in A\}} 
\frac{\Gamma(Q)}{\e+\Gamma(Q)}
d\mu=0.
\end{equation}

\bigskip

{\bf Step 3}. Observe that $\Gamma(Q)$ is a polynomial of degree at most $2d$, see indeed (\ref{gamma1}), (\ref{gamma2}) or (\ref{gamma3}). We deduce from the Carbery-Wright inequality (\ref{cw-ineq}), together with the Poincar\'e inequality (\ref{poincare}),
that $\Gamma(Q)$ is strictly positive almost everywhere. Thus, by dominated convergence $(\e\to 0)$ in (\ref{bla}) we finally get that $\mu\{Q\in A\}=P(Q(X_1,\ldots,X_{m})\in A)=0$.
\qed

\section{Proof of Theorem \ref{Main}}

We are now in a position to show Theorem \ref{Main}.
We will split its proof in several steps.\\

{\bf Step 1}. For any $p\in [1,\infty)$ we shall prove that
\begin{equation}\label{claim1}
\sup_n \int |Q_n|^pd\mu_{m}<\infty.
\end{equation}
(Let us recall our convention about $m$ from the beginning of Section 2.)
Indeed, using (for instance) Propositions 3.11, 3.12 and 3.16 of \cite{MOO} (namely, a hypercontractivity property), one
first observes that, for any $p\in[2,\infty)$,
there exists a constant
$c_p>0$ such that, for  all $n$,
\begin{equation}\label{hypercontractivity}
\int |Q_n|^pd\mu_{m}\leq c_p\, \left(\int Q_n^2d\mu_{m}\right)^{p/2}.
\end{equation}
(This is  for obtaining (\ref{hypercontractivity}) that
we need $Q_n$ to be {\it multilinear}.)
On the other hand, one can write
\begin{eqnarray*}
\int Q_n^2d\mu_{m}&=&\int Q_n^2\,{\bf 1}_{\{Q_n^2\geq\frac12\int Q_n^2d\mu_{m}\}}d\mu_{m}+\int Q_n^2\,{\bf 1}_{\{Q_n^2<\frac12\int Q_n^2d\mu_{m}\}}d\mu_{m}\\
&\leq& \sqrt{\int Q_n^4d\mu_{m}}\sqrt{\mu_{m}\left\{x:Q_n(x)^2\geq\frac12\int Q_n^2d\mu_{m}\right\}}+\frac12 \int Q_n^2d\mu_{m},
\end{eqnarray*}
so that, using (\ref{hypercontractivity}) with $p=4$,
\[
\mu_{m}\left\{x:\,Q_n(x)^2\geq\frac12\int Q_n^2d\mu_{m}\right\}
\geq \frac{\left(\int Q_n^2d\mu_{m}\right)^2}{4\int Q_n^4d\mu_{m}}\geq \frac1{4c_4}.
\]
But $\{Q_n\}_{n\geq 1}$ is tight as $\{F_n\}_{n\geq 1}$ converges in law. As a result, there exists $M>0$ such that,  for all $n$,
\[
\mu_{m}\left\{x:\,Q_n(x)^2\geq M\right\}
<\frac1{4c_4}.
\]
We deduce that $\int Q_n^2d\mu_{m}\leq 2M$ which, together with (\ref{hypercontractivity}), leads to the claim
(\ref{claim1}).

\bigskip
\bigskip

{\bf Step 2}. We shall prove the existence of a constant $c>0$ such that, for any $u>0$ and any $n\in\N^*$,
\begin{equation}\label{claim2}
\mu_{m}\left\{x:\,\Gamma(Q_n)\leq u\right\}\leq c\, \frac{u^{\frac{1}{2d}}}{{\rm Var}_{\mu_{m}}(Q_n)^{\frac{1}{2d}}}.
\end{equation}
Observe first that $\Gamma(Q_n)$ is a polynomial of degree at most $2d$, see indeed (\ref{gamma1}), (\ref{gamma2}) or (\ref{gamma3}). On the other hand, since $X$ has a log-concave density, the probability $\mu_{m}$ is absolutely continuous
with a log-concave density as well.
As a consequence, Carbery-Wright inequality (\ref{cw-ineq}) applies and yields the existence of a constant $c>0$ such that
\[
\mu_{m}\left\{x:\,\Gamma(Q_n)\leq u\right\}\leq c\,u^{\frac{1}{2d}}\,
\left(
\int \Gamma(Q_n)d\mu_{m}
\right)
^{-\frac{1}{2d}}.
\]
To get the claim  (\ref{claim2}), it remains one to apply the Poincar\'e inequality (\ref{poincare}).

\bigskip

{\bf Step 3}. We shall prove the existence of $n_0 \in\N^*$ and $\kappa>0$ such
that, for any $\e>0$,
\begin{equation}\label{claim3}
\sup_{n\geq n_0} \int \frac{\e}{\Gamma(Q_n)+\e}d\mu_{m}\leq \kappa\,\e^{\frac{1}{2d+1}}.
\end{equation}
Indeed, thanks to the result shown in Step 2 one can write
\begin{eqnarray*}
\int \frac{\e}{\Gamma(Q_n)+\e}d\mu_{m}&\leq& \frac{\e}{u} + \mu_{m}\left\{x:\,\Gamma(Q_n)\leq u\right\}\\
&\leq& \frac{\e}{u} +c\, \frac{u^{\frac{1}{2d}}}{{\rm Var_{\mu_{m}}}(Q_n)^{\frac{1}{2d}}}.
\end{eqnarray*}
But, by Step 1 and since $\mu_{m}\circ Q_n^{-1}$ converges to some probability measure $\eta$, one has that ${\rm Var}_{\mu_{m}}(Q_n)$ converges to the variance of $\eta$ as $n\to\infty$. Moreover this variance is strictly positive
by assumption. We deduce the existence of $n_0 \in\N^*$ and $\delta>0$ such that
\[
\sup_{n\geq n_0} \int \frac{\e}{\Gamma(Q_n)+\e}d\mu_{m}
\leq \frac{\e}{u} +\delta\, u^{\frac{1}{2d}}.
\]
Choosing $u=\e^{\frac{2d}{2d+1}}$ leads to the desired conclusion (\ref{claim3}).

\bigskip

{\bf Step 4}. Let $m'$ be shorthand for $m_{n'}$ and 
recall the Fortet-Mourier distance (\ref{fm}) as well as the total variation distance (\ref{tv2}) from the Introduction.
We shall prove that, for any $n,n'\geq n_0$ (with $n_0$ and $\kappa$ given by Step 3), any $0< \alpha\leq 1$ and any $\e>0$, 
\begin{equation}\label{claim4}
d_{TV}(F_n,F_{n'})\leq \frac{1}{\alpha}d_{FM}(F_n,F_{n'})+4\kappa\,\e^{\frac{1}{2d+1}}
+2\sqrt{\frac{2}\pi}\,\frac{\alpha}{\e}\,\sup_{n\geq n_0}
\left(
\int \Gamma(\Gamma(Q_n))d\mu_{m} + \int \big|\mathcal{L}Q_n\big|d\mu_{m}
\right).
\end{equation}
Indeed, set $p_{\alpha}(x)=\frac{1}{\alpha\sqrt{2\pi}}e^{-\frac{x^2}{2\alpha^2}}$, $x\in\R$, $0<\alpha\leq 1$, and
let $g\in C^\infty_c$ be bounded by 1.
It is immediately checked that
\begin{equation}\label{fc1}
\|g*p_\alpha\|_\infty\leq 1\leq\frac1\alpha
\quad\mbox{and}\quad
\|(g*p_\alpha)'\|_\infty
\leq\frac1\alpha.
\end{equation}
Let $n,n'\geq n_0$ be given integers. 
Using Step 3 and (\ref{fc1})
we can write
\begin{eqnarray*}
&&\left|\int g\, d(\mu_{m}\circ Q_n^{-1})-\int g\,d(\mu_{m'}\circ Q_{n'}^{-1})\right|=
\left|\int g\circ Q_n\,d\mu_{m}-\int g\circ Q_{n'}\,d\mu_{m'}\right|\\
&\leq&\left|\int(g*p_\alpha)\circ Q_nd\mu_{m}-\int (g*p_\alpha)\circ Q_{n'}d\mu_{m'}\right|\\
&&
+\left|\int(g-g*p_\alpha)\circ Q_n\times 
\left(\frac{\Gamma(Q_{n})}{\Gamma(Q_n)+\varepsilon}+\frac{\e}{\Gamma(Q_n)+\varepsilon}\right)d\mu_{m}\right|\\
&&
+\left|\int(g-g*p_\alpha)\circ Q_{n'}\times
\left(\frac{\Gamma(Q_{n'})}{\Gamma(Q_{n'})+\varepsilon}+\frac{\e}{\Gamma(Q_{n'})+\varepsilon}\right)d\mu_{m'}\right|\\
&\leq & \frac{1}{\alpha}d_{FM}(F_n,F_{n'})
+2\int \frac{\varepsilon}{\Gamma(Q_{n})+\varepsilon}d\mu_{m}
+2\int \frac{\varepsilon}{\Gamma(Q_{n'})+\varepsilon}d\mu_{m'}
\\
&&
+\left|\int(g-g*p_\alpha)\circ Q_{n}\times
\frac{\Gamma(Q_{n})}{\Gamma(Q_{n})+\varepsilon}d\mu_{m}\right|
+\left|\int(g-g*p_\alpha)\circ Q_{n'}\times
\frac{\Gamma(Q_{n'})}{\Gamma(Q_{n'})+\varepsilon}d\mu_{m'}\right|
\\
&\leq & \frac{1}{\alpha}d_{FM}(F_n,F_{n'})+4\kappa\,\e^{\frac{1}{2d+1}}
+2\sup_{n\geq n_0}\left|\int(g-g*p_\alpha)\circ Q_{n}\times
\frac{\Gamma(Q_{n})}{\Gamma(Q_{n})+\varepsilon}d\mu_{m}\right|.
\end{eqnarray*}
Now, set $\Psi(x)=\int_{-\infty}^{x} g(s) ds$ and let us apply (\ref{diffu}). We obtain
\begin{eqnarray}
\notag
&&\left|\int(g-g*p_\alpha)\circ Q_{n}\times
\frac{\Gamma(Q_{n})}{\Gamma(Q_{n})+\varepsilon}d\mu_{m}\right|\\
&=&\left|\int \frac{1}{\Gamma(Q_{n})+\varepsilon}\,\Gamma\big((\Psi-\Psi*p_{\alpha})\circ Q_n,Q_n\big)d\mu_{m}\right|\notag\\
\notag
&=&\left|\int (\Psi-\Psi*p_{\alpha})\circ Q_n\times \left(\Gamma\big(Q_n,\frac{1}{\Gamma(Q_{n})+\varepsilon}\big)+\frac{\mathcal{L}Q_n}{\Gamma(Q_{n})+\varepsilon}\right)
d\mu_{m}\right|\\\notag
&=&\left|\int(\Psi-\Psi*p_{\alpha})\circ Q_n\times \left(-\frac{\Gamma(Q_n,\Gamma(Q_n))}
{(\Gamma(Q_{n})+\varepsilon)^2}+\frac{\mathcal{L}Q_n}{\Gamma(Q_{n})+\varepsilon}\right)d\mu_{m}\right|\\
&\leq&\frac1\e\,\int|(\Psi-\Psi*p_{\alpha})\circ Q_n|\times\big(\Gamma(\Gamma(Q_n))+\big|\mathcal{L}Q_n\big|\big)d\mu_{m_{n}}.\label{FactA}
\end{eqnarray}
On the other hand, we have
\begin{eqnarray}
\notag
\left|\Psi(x)-\Psi*p_{\alpha}(x)\right|&=&\left|\int_\R p_\alpha(y)\left(\int_{-\infty}^x \left(g(u)-g(u-y)\right)du\right)dy\right|\\\notag
&\leq&\int_\R p_\alpha(y)\left|\int_{-\infty}^x g(u)du-\int_{-\infty}^xg(u-y)du\right|dy\\
&\leq&\int_\R p_\alpha(y)\left|\int_{x-y}^x g(u)du\right| dy
\leq \int_\R p_{\alpha}(y)\left|y\right|dy \leq \sqrt{\frac{2}{\pi}}\alpha.\label{FactB}
\end{eqnarray}
The desired conclusion (\ref{claim4}) now follows easily.

\bigskip
\bigskip

{\bf Step 5}. We shall prove that
\begin{equation}\label{claim5}
\sup_{n\geq n_0}
\left(
\int \Gamma(\Gamma(Q_n))d\mu_{m} + \int \big|\mathcal{L}Q_n\big|d\mu_{m}
\right)<\infty.
\end{equation}
First, relying on the results of Section \ref{preli} we have that
\[
Q_n \in \bigoplus_{\alpha\leq \lambda_{2d}}\text{Ker}(\mathcal{L}+\alpha I).
\]
Since $\mathcal{L}$ is a bounded operator on the space $\bigoplus_{\alpha\leq \lambda_{2d}}\text{Ker}(\mathcal{L}+\alpha I)$ and $Q_n$ is bounded in the space $L^2(\mu_m)$, we deduce immediately that $\sup_n\int \mathcal{L}(Q_n)^2d\mu_{m}<\infty$,
implying in turn that
\[
\sup_n\int |\mathcal{L}(Q_n)|d\mu_{m}<\infty.
\]
Besides, one has $\Gamma=\frac{1}{2}(\mathcal{L}+2\lambda I)$
on $\text{Ker}(\mathcal{L}+\lambda I)$ and one deduces for the same reason as above that
\[
\sup_{n}\int\Gamma(\Gamma(Q_n))d\mu_{m}<\infty.
\]
The proof of (\ref{claim5}) is complete.
\bigskip
\bigskip

{\bf Step 6: conclusion}. 
The Fortet-Mourier distance $d_{FM}$ metrizing the convergence in distribution, our assumption ensures that $d_{FM}(F_n,F_{n'})\to 0$ as $n,n'\to\infty$.
Therefore, combining (\ref{claim5}) with (\ref{claim4}), letting $n,n'\to\infty$, then $\alpha\to 0$ and
then $\e\to 0$, we conclude that $\lim_{n,n'\to \infty} d_{TV}(F_n,F_{n'})=0$, meaning that
$F_n$ is a Cauchy sequence in the total variation topology.
But the space of bounded measures is complete for the total variation distance, so the distribution of $F_n$ must converge to some distribution, say $\eta$, in the total variation distance. 
Of course, $\eta$ must coincide with the law of $F_\infty$.
Moreover, let $A$ be a Borel set of Lebesgue measure zero.
By Proposition \ref{abso}, we have $P(F_n\in A)=0$ when $n$ is large enough.
Since $d_{TV}(F_n,F_\infty)\to 0$ as $n\to\infty$, we deduce that $P(F_\infty\in A)=0$ as well, thus
proving that the law of $F_\infty$ is absolutely continuous with respect to the Lebesgue measure by
the Radon-Nikodym theorem.
The proof of Theorem \ref{Main} is now complete.
\qed

\end{document}